\documentclass[a4paper,11pt]{article}
\usepackage{amssymb}
\usepackage{amsmath}
\usepackage{tikz-cd}

\usepackage{amsthm}
\usepackage{amsfonts}
\usepackage{mathrsfs}
\usepackage{afterpage}
\usepackage{float}
\usepackage{multirow}
\usepackage{graphicx}
\usepackage{color}
\usepackage[margin=1.2in]{geometry}
\newcommand*{\SHom}{\mathscr{H}\kern -.5pt om}
\newcommand*{\SExt}{\mathscr{E}\kern -.5pt xt}

\begin{document}
\title{\textbf{Stability Conditions and Exceptional Objects in Triangulated Categories}}
\date{}
\author{Zihong Chen}
\maketitle
\theoremstyle{definition}
\newtheorem{mydef}{Definition}[section]
\newtheorem{rmk}{Remark}[section]
\theoremstyle{plain}
\newtheorem{cor}{Corollary}[section]
\newtheorem{lemma}{Lemma}[section]
\newtheorem{thm}{Theorem}[section]
\newtheorem{prop}{Proposition}[section]

\begin{abstract}
The goal of this paper is to study the subspace of stability condition $\Sigma_{\mathcal{E}}\subset \mathrm{Stab}(X)$ associated to an exceptional collection $\mathcal{E}$ on a projective variety $X$. Following Macr\`{i}'s approach, we show a certain correspondence between the homotopy class of continuous loops in $\Sigma_{\mathcal{E}}$ and words of the braid group. In particular, we prove that in the case $X=\mathbb{P}^3$ and $\mathcal{E}=\{\mathcal{O},\mathcal{O}(1),\mathcal{O}(2),\mathcal{O}(3)\}$, the space $\Sigma_{\mathcal{E}}$ is a connected and simply connected 4-dimensional manifold.
\end{abstract}
\section{Introduction}
T.Bridgeland introduced the notion of stability condition on a triangulated category in \cite{Bri07}. The motivation came from Douglas's work on $\Pi$-stability in string theory (\cite{Dou02}). Bridgeland's stability condition also generalizes the $\mu$-stability for coherent sheaves on projective varieties. \par\indent
One main result of Bridgeland is that the set of all stability conditions $\mathrm{Stab}(\mathcal{D})$ form a topological space (\cite[Theorem 1.2]{Bri07}). Provided that certain technical conditions are satisfied (local finiteness), $\mathrm{Stab}(\mathcal{D})$ is furthermore a smooth manifold.\par\indent
A few explicit computations of $\mathrm{Stab}(X):=\mathrm{Stab}(D^b(\mathrm{Coh}(X)))$ for a smooth projective variety $X$ have been done. Bridgeland showed that when $C$ is an elliptic curve, $\mathrm{Stab}(C)$ is isomorphic to $\widetilde{\mathrm{GL}^+(2,\mathbb{R})}$ through a free and transitive action of the latter on the former (\cite[Section 9]{Bri07}). E.Macr\`{i} generalized this result to all smooth projective curves with genus greater or equal to $1$ (\cite{Mac07}). The case of $X=\mathbb{P}^1$ was computed by S.Okada, who showed that $\mathrm{Stab}(\mathbb{P}^1)\cong \mathbb{C}^2$. The case of several special surfaces have also been studied in \cite{Bri08},\cite{BM11}. In general, such explicit results are hard to obtain, but certain common properties of the topological space $\mathrm{Stab}(\mathcal{D})$ have been found among various examples. For instance, it has been conjectured that each connected component $\Sigma\subset \mathrm{Stab}(\mathcal{D})$ is simply connected or even contractible.\par\indent
In \cite{Mac07}, Macr\`{i} studied stability conditions generated by a finite complete exceptional collection of objects. To each complete exceptional collection $\mathcal{E}=\{E_0,\cdots,E_n\}$, he associates an open subspace of stability conditions for which $E_0,\cdots,E_n$ are stable, denoted by $\Theta_{\mathcal{E}}$. He showed that each $\Theta_{\mathcal{E}}$ is a connected and simply connected open submanifold of $\mathrm{Stab}(\mathcal{D})$ with maximal dimension. Furthermore, the union of $\Theta_{\mathcal{F}}$ as $\mathcal{F}$ ranges through all iterated mutations of $\mathcal{E}$, denoted by $\Sigma_{\mathcal{E}}$, is again connected. In the special case when $\mathcal{D}=D^b(\mathrm{Coh}(\mathbb{P}^1))$ and $\mathcal{D}=D^b(\mathrm{Coh}(\mathbb{P}^2))$, Macr\`{i} showed that $\Sigma_{\mathcal{E}}$ is in fact simply connected for each $\mathcal{E}$. \\\par\indent
In this paper, we study the simply connectedness of $\Sigma_{\mathcal{E}}$ in a more general context. In particular, we prove the following proposition, which generalizes Macr\`{i}'s result to an arbitrary triangulated category equipped with an exceptional collection satisfying certain conditions (denoted by $\dag$, see Section 4).
\begin{prop}
Fix a triangulated category $\mathcal{D}$. Let $\mathcal{E}=\{E_0,\cdots,E_n\}$ be an exceptional collection satisfying $\dag$. Let $\gamma: [0,1]\rightarrow \Sigma_{\mathcal{E}}$ be a continuous loop with $\gamma(0)=\gamma(1)\in \Theta_{\mathcal{E}}$. Then, up to replacing $\gamma$ by a homotopic path, there exists $l=l_s\cdots l_1$ with $l_i\in \{\mathcal{L}_0,\cdots,\mathcal{L}_{n-1},\mathcal{R}_0,\cdots,\mathcal{R}_{n-1}\}$ for all $i$, such that $l\mathcal{E}=\mathcal{E}$, and real numbers $0=a_0<a_1<\cdots<a_s<a_{s+1}=1$ such that $\gamma([a_k,a_{k+1}))\subset \Theta_{l_k\cdots l_1\mathcal{E}}$ for all $k=0,1,\cdots,s$.
\end{prop}
The motivation for studying exceptional collections satisfying $\dag$ comes from the study of coherent sheaves over projective varieties. A well known such example is the collection $\mathcal{E}=\{\mathcal{O},\mathcal{O}(1),\cdots,\mathcal{O}(n)\}$ in $D^b(\mathrm{Coh}(\mathbb{P}^n))$; see \cite{Bon90}.\par\indent
Using Proposition 1.1, and explicit computations regarding the braid action on exceptional collections on $\mathbb{P}^3$, we prove the following main theorem.
\begin{thm}
Let $\mathcal{E}$ be the exceptional collection $\{\mathcal{O},\mathcal{O}(1),\mathcal{O}(2),\mathcal{O}(3)\}$ on $\mathbb{P}^3$. Then, the subspace $\Sigma_{\mathcal{E}}$ of $\mathrm{Stab}(\mathbb{P}^3)$ is simply connected.
\end{thm}
The organization of this paper is as follows: In Section 2 we briefly review some aspects of Bridgeland stability conditions. In Sections 3 we summarize some of Macr\`{i}'s concepts on stability conditions generated by a finite exceptional collection. We build up the proof of Proposition 1.1 in Section 4, and in Section 5 we study the case of $\mathbb{P}^3$ and prove our main theorem. In the appendix, we review some basic facts about homological algebra and algebraic geometry.\\\\
\textbf{Acknowledgements.} I would like to thank Hiro Lee Tanaka for introducing me to this topic and his guidance throughout my research. I would also like to thank Emanuele Macr\`{i} for his invaluable advice and comments.

\section{Stability Conditions on Triangulated Categories}
This section is a brief summary of \cite{Bri07} and serve as a review for some basic concepts leading to the definition of Bridgeland Stability. In the discussion below, we assume that all triangulated categories are small and Hom-finite over a fixed field $K$ (i.e. $\mathrm{Hom}(A,B)$ is a finite dimensional vector space for all $A,B\in \mathcal{D}$). Let $K(\mathcal{D})$ denote the Grothedieck group of $\mathcal{D}$. \par\indent
First, we recall the definition of a $t$-structure.
\theoremstyle{definition}
\begin{mydef}
A \emph{$t$-structure} on a triangulated category $\mathcal{D}$ is the data of a pair of full subcategories $(D^{\leq 0}, D^{\geq 0})$ satisfying the following
conditions:\\
(TS1) If we denote $\mathcal{D}^{\leq n}=\mathcal{D}^{\leq 0}[-n]$, and $\mathcal{D}^{\geq n}=\mathcal{D}^{\geq 0}[-n]$, then $\mathcal{D}^{\leq 0}\subset \mathcal{D}^{\leq 1}$ and $\mathcal{D}^{\geq 0}\supset\mathcal{D}^{\geq 1}$.\\
(TS2) For any $x\in \mathcal{D}^{\leq 0}, y\in \mathcal{D}^{\geq 1}$, we have $\mathrm{Hom}(x,y)=0$.\\
(TS3) For any $x\in \mathcal{D}$, there exists an exact triangle $x_{\leq 0} \longrightarrow x\longrightarrow x_{\geq 1}$, with $x_{\leq0}\in \mathcal{D}^{\leq 0}$ and $x_{\geq 1}\in \mathcal{D}^{\geq
1}$.
\end{mydef}
We define the \emph{heart} of this $t$-structure to be $\mathcal{A}=D^{\leq 0}\bigcap D^{\geq 0}$, which turns out to always be an abelian category. We say that the $t$-structure is \emph{bounded} if $\mathcal{D}=\bigcup_{i,j\in \mathbb{Z}} \mathcal{D}^{\leq i}\bigcap \mathcal{D}^{\geq j}$. From now on, unless mentioned otherwise, we always assume that a $t$-structure is bounded.\par\indent
A bounded $t$-structure is uniquely determined by its heart, and therefore we can interchange these two concepts. The following lemma tells us that the heart of a bounded $t$-structure generalizes the concept `filtration by cohomology'.
\begin{lemma}
Let $\mathcal{A}\subset D$ be a full additive subcategory of a triangulated category $\mathcal{D}$. Then $\mathcal{A}$ is the heart of a bounded $t$-structure on $\mathcal{D}$ if and only if the following two conditions hold:\\
(a) for integers $k_1>k_2$, we have $\mathrm{Hom}(A[k_1], B[k_2])=0$ for all $A,B\in \mathcal{A}$.\\
(b) for every nonzero object $E\in D$ there is a finite sequence of integers
$$k_1>k_2>\cdots>k_n$$
and a collection of triangles
\begin{center}
\begin{tikzcd}[row sep=1.2cm, column sep=0.8cm]
0=E_0\arrow[rr] & &E_1\arrow[rr]\arrow[dl]& &E_2\arrow[r]\arrow[dl]&\cdots\arrow[r]&E_{n-1}\arrow[rr]& &E_{n}=E\arrow[dl]\\
& A_1\arrow[ul]& & A_2\arrow[ul]& & & & A_n\arrow[ul] &
\end{tikzcd}
\end{center}
with $A_k\in \mathcal{A}[k_j]$.
\end{lemma}
Now, we give the definition of a stability function on an abelian category (see \cite{Rud97}), which is historically prior to Bridgeland's notion of stability conditions on triangulated categories.
\begin{mydef}
A \emph{stability function} on an abelian category $\mathcal{A}$ is a group homomorphism $Z: K(\mathcal{A})\rightarrow \mathbb{C}$ such that for all $0\neq E\in \mathcal{A}$, the complex number $Z(E)$ lies in the strict upper half plane
$$H=\{re^{i\pi \phi}: r>0, \phi\in(0,1]\}.$$
\end{mydef}
Given a stability function $Z$, the \emph{phase} of a nonzero object $E$ is defined as $\phi(E):=(1/\pi)\mathrm{arg}(Z(E))\in (0,1]$. $E$ is called \emph{semistable} (\emph{stable}) if every $0\neq A\hookrightarrow E$ satisfies $\phi(A)\leq \phi(E)$ ($\phi(A)<\phi(E)$). Central to the study of stability functions is the notion of a \emph{Harder-Narasimhan filtration}.
\begin{mydef}
Given a stability function $Z: K(\mathcal{A})\rightarrow \mathbb{C}$, a \emph{Harder-Narasimhan filtration} of a nonzero objects $0\neq E\in \mathcal{A}$ is a finite chain of subobjects
$$0=E_0\hookrightarrow E_1\hookrightarrow \cdots\hookrightarrow E_{n-1}\hookrightarrow E_n=E$$
such that $F_i=E_i/E_{i-1}$ are semistable and
$$\phi(F_1)>\phi(F_2)>\cdots>\phi(F_n).$$
If such a filtration exists, the stability function is said to have the \emph{Harder-Narasimhan property}.
\end{mydef}
As an example, let $C$ be a smooth projective curve. Then, the stability function $Z$ on $\mathrm{Coh}(C)$ defined by $Z(E)=-\mathrm{deg}(E)+i\,\mathrm{rank}(E)$ satisfies the Harder-Narasimhan property. \\\\
\emph{Remark}. Note that both Lemma 2.1 and Definition 2.3 involve certain kind of filtration by triangles. However, while the former filtration ranges across different hearts, the latter lies within a single abelian category. In some sense, a stability condition on a triangulated category is a combination of the two, an intuition that is made precise by Proposition 2.1.
\begin{mydef}
A \emph{stability condition} $\sigma=(Z,\mathcal{P})$ on a triangulated category $\mathcal{D}$ consists of a group homomorphism $Z: K(\mathcal{D})\rightarrow \mathbb{C}$ called the \emph{central charge}, and full additive subcategories $\mathcal{P}(\phi)$ for each $\phi\in\mathbb{R}$, satisfying the following conditions:\\
(a) if $0\neq E\in \mathcal{P}(\phi)$ then $Z(E)\in \mathbb{R}_{>0}e^{i\pi\phi}$,\\
(b) for all $\phi\in \mathbb{R}$, $\mathcal{P}(\phi+1)=\mathcal{P}(\phi)[1]$,\\
(c) if $\phi_1>\phi_2$ and $A_j\in \mathcal{P}(\phi_j)$ then $\mathrm{Hom}(A_1,A_2)=0$,\\
(d) for each nonzero $E\in\mathcal{D}$, there exists a finite sequence of real numbers $\phi_1>\phi_2>\cdots>\phi_n$ and a collection of triangles
\begin{center}
\begin{tikzcd}[row sep=1.2cm, column sep=0.8cm]
0=E_0\arrow[rr] & &E_1\arrow[rr]\arrow[dl]& &E_2\arrow[r]\arrow[dl]&\cdots\arrow[r]&E_{n-1}\arrow[rr]& &E_{n}=E\arrow[dl]\\
& A_1\arrow[ul]& & A_2\arrow[ul]& & & & A_n\arrow[ul] &
\end{tikzcd}
\end{center}
with $A_j\in\mathcal{P}(\phi_j)$ for all $j$.
\end{mydef}
In the above definition, the filtration in (4) is also called a Harder-Narasimhan filtration, which is unique up to isomorphism. Hence, we may define $\phi^+_{\sigma}(E):=\phi_1,\,\phi^-_{\sigma}(E):=\phi_n$ and $m_{\sigma}(E)=\sum_j|Z(A_j)|$. The nonzero objects of $\mathcal{P}(\phi)$ are called \emph{semistable} in $\sigma$ of phase $\phi$; the simple objects of $\mathcal{P}(\phi)$ are called \emph{stable}. In fact, each $\mathcal{P}(\phi)$ is an abelian category (\cite[Lemma 5.2]{Bri07}).\par\indent
For an interval $I$, let $\mathcal{P}(I)$ denote the extension-closed subcategory generated by $\mathcal{P}(\phi)$ for $\phi\in I$. We call the abelian category $\mathcal{P}((0,1])$ the \emph{heart} of the stability condition; in fact, it is the heart of the $t$-structure $(\mathcal{P}(>0), \mathcal{P}(\leq 1))$.
\begin{prop}
To give a stability condition on a triangulated category $\mathcal{D}$ is equivalent to giving a bounded $t$-structure on $\mathcal{D}$ and a stability function on its heart with the Harder-Narasimhan property.
\end{prop}
This proposition implies, for instance, that given a smooth projective curve $C$, the stability function on $\mathrm{Coh}(C)$ given by $Z(E)=-\mathrm{deg}(E)+i\,\mathrm{rank}(E)$ induces a stability condition on $D^b(\mathrm{Coh}(C))$.\par\indent
A stability condition $(Z,\mathcal{P})$ is called \emph{locally finite} if for each $\phi\in \mathbb{R}$, there exists $\epsilon>0$ such that the quasi-abelian category $\mathcal{P}(\phi-\epsilon,\phi+\epsilon)$ is of finite length. In particular, this implies that $\mathcal{P}(\phi)$ is of finite length, and hence every semistable object of phase $\phi$ has a finite Jordan-H\"{o}lder filtration with stable factors of the same phase. We denote the set of all locally finite stability conditions on $\mathcal{D}$ by $\mathrm{Stab}(\mathcal{D})$.\par\indent
One of the most important feature of $\mathrm{Stab}(\mathcal{D})$ is its natural topology defined as follows. For $\sigma_1,\sigma_2\in\mathrm{Stab}(\mathcal{D})$, the function
$$d(\sigma_1,\sigma_2)=\sup_{0\neq E\in\mathcal{D}}\Big\{|\phi^-_{\sigma_1}(E)-\phi^-_{\sigma_2}(E)|, |\phi^+_{\sigma_1}(E)-\phi^+_{\sigma_2}(E)|, |\log\frac{m_{\sigma_1}(E)}{m_{\sigma_2}(E)}|\Big\}$$
defines a generalized metric on $\mathrm{Stab}(\mathcal{D})$. \par\indent
Given a stability condition $\sigma=(Z,\mathcal{P})$, the natural projection $(Z,\mathcal{P})\mapsto Z$ induces a continuous map from $\mathrm{Stab}(\mathcal{D})$ to $\mathrm{Hom}_{\mathbb{Z}}(K(\mathcal{D},\mathbb{C})$. Bridgeland proved that this map in fact a local homeomorphism.
\begin{thm}
Let $\mathcal{D}$ be a triangulated category. For each connected component $\Sigma\subset \mathrm{Stab}(\mathcal{D})$ there is a linear subspace $V(\Sigma)\subset \mathrm{Hom}_{\mathbb{Z}}(K(\mathcal{D}),\mathbb{C})$, with a well defined linear topology, and a local homeomorphism $\mathcal{Z}:\Sigma\rightarrow V(\Sigma)$ sending $(Z,\mathcal{P})$ to its central charge $Z$.
\end{thm}
Finally, we remark that $\mathrm{Stab}(\mathrm{D})$ carries a right action by $\widetilde{\mathrm{GL}^+(2,\mathbb{R})}$, the universal cover of $\mathrm{GL}^+(2,\mathbb{R})$, and a left action by $\mathrm{Aut}(\mathcal{D})$, the group of exact autoequivalences of $\mathcal{D}$ (\cite[Lemma 8.2]{Bri07}). \par\indent
First note that we can explicitly write $\widetilde{\mathrm{GL}^+(2,\mathbb{R})}=$
$$\{(T,f): T\in \mathrm{GL}^+(2,\mathbb{R}), f: \mathbb{R}\rightarrow \mathbb{R}\;\mathrm{increasing\;with}\;f(\phi+1)=f(\phi)+1, \mathrm{and}\;Te^{i\pi\phi}\in \mathbb{R}_{>0}e^{i\pi f(\phi)}\}.$$
Given $(T,f)\in \widetilde{\mathrm{GL}^+(2,\mathbb{R})}$ and $(Z,\mathcal{P})\in \mathrm{Stab}(\mathcal{D})$, we define the action by $(T,f)\cdot (Z,\mathcal{P})=(T^{-1}\circ Z, \mathcal{P}\circ f)$. In essence, an action of $(T,f)$ is a relabeling of the phase of $(Z,\mathcal{P})$ (with some rescaling), but the set of semistable (stable) objects are left unchanged. Therefore, it is often convenient to identify two stability conditions up to the $\widetilde{\mathrm{GL}^+(2,\mathbb{R})}$ action. \par\indent
Finally, given $\Psi\in \mathrm{Aut}(\mathcal{D})$, let $\psi$ denote the induced map on $(K(\mathcal{D})$. Then we define the action by $\Psi\cdot (Z,\mathcal{P})=(Z\circ \psi, \Psi\circ \mathcal{P})$. It is clear that this action commutes with the action of $\widetilde{\mathrm{GL}^+(2,\mathbb{R})}$.

\section{Some Properties of Exceptional Objects}
In this section, we review the basics of exceptional objects following \cite[Section 2]{Bon90}. We also discuss Macr\`{i}'s approach to stability conditions via exceptional collections in \cite{Mac04},\cite{Mac07}. \par\indent
As before, let $\mathcal{D}$ be a small and $\mathrm{Hom}$-finite triangulated category linear over some field $K$. For $A,B\in\mathcal{D}$, we define their $\mathrm{Hom}$ complex to be
$$\mathrm{Hom}^{\bullet}(A,B):=\bigoplus_{k\in\mathbb{Z}}\mathrm{Hom}^k(A,B)[-k],$$
where $\mathrm{Hom}^k(A,B):=\mathrm{Ext}^k(A,B):=\mathrm{Hom}(A,B[k])$. \par\indent
\begin{mydef}
(i) An object $E\in\mathcal{D}$ is called an \emph{exceptional object} if $\mathrm{Hom}^i(E,E)=0$ when $i\neq 0$ and $\mathrm{Hom}(E,E)=K$.\\
(ii) An ordered collection of exceptional objects $\{E_0,\cdots,E_n\}$ is called an \emph{exceptional collection} in $\mathcal{D}$ if for all $j>i$, $\mathrm{Hom}^{\bullet}(E_j,E_i)=0$. An exceptional collection consisting of two elements is called an \emph{exceptional pair}.
\end{mydef}
\begin{mydef}
Let $\mathcal{E}=\{E_0,\cdots,E_n\}$ be an exceptional collection. $\mathcal{E}$ is called\\
$\bullet$ strong, if $\mathrm{Hom}^k(E_i,E_j)=0$ for all $i,j$, with $k\neq 0$;\\
$\bullet$ Ext, if $\mathrm{Hom}^{\leq 0}(E_i,E_j)=0$ for all $i\neq j$;\\
$\bullet$ complete, if $\mathcal{E}$ generates $\mathcal{D}$ by shifts and extensions.
\end{mydef}
\begin{mydef}
Let $\{E,F\}$ be an exceptional pair. We define the \emph{left mutation} $\mathcal{L}_EF$ and the \emph{right mutation} $\mathcal{R}_FE$ with the aid of distinguished triangles in $\mathcal{D}$:
$$\mathcal{L}_EF\rightarrow \mathrm{Hom}^{\bullet}(E,F)\otimes E\rightarrow F,$$
$$E\rightarrow \mathrm{Hom}^{\bullet}(E,F)^*\otimes F\rightarrow \mathcal{R}_FE,$$
where $V[k]\otimes E$ is defined as $E[k]^{\oplus \dim V}$. Note that under duality of vector spaces the grading changes sign.
\end{mydef}
A \emph{mutation} of an exceptional collection $\mathcal{E}=\{E_0,\cdots,E_n\}$ is defined as a mutation of a pair in this collection:
$$\mathcal{R}_i\mathcal{E}=\{E_0,\cdots,E_{i-1},E_{i+1},\mathcal{R}_{E_{i+1}}E_i,E_{i+2},\cdots,E_n\},$$
$$\mathcal{L}_i\mathcal{E}=\{E_0,\cdots,E_{i-1},\mathcal{L}_{E_i}E_{i+1},E_i,E_{i+2},\cdots,E_n\},$$
for $i=0,1,\cdots,n-1$. By the following proposition (see \cite[Section 2]{Bon90}), a mutation of an exceptional collection is still exceptional, and thus we may define mutations on the mutated collection. Composition of mutations constructed in this way is called an \emph{iterated mutation}.
\begin{prop}
(i) A mutation of an exceptional collection is an exceptional collection.\\
(ii) A mutation of a complete exceptional collection is complete exceptional.\\
(iii) The following relations hold:
$$\mathcal{R}_i\mathcal{L}_i=\mathcal{L}_i\mathcal{R}_i=1\qquad \mathcal{R}_i\mathcal{R}_{i+1}\mathcal{R}_i=\mathcal{R}_{i+1}\mathcal{R}_i\mathcal{R}_{i+1}\qquad \mathcal{L}_i\mathcal{L}_{i+1}\mathcal{L}_i=\mathcal{L}_{i+1}\mathcal{L}_i\mathcal{L}_{i+1}.$$
\end{prop}
\noindent\emph{Remark}. Recall that the $(n+1)$-th Artin braid group $A_{n+1}$ can be defined via the presentation $$A_{n+1}=\langle\sigma_0,\cdots,\sigma_{n-1}\,|\,\sigma_i\sigma_{i+1}\sigma_i=\sigma_{i+1}\sigma_{i}\sigma_{i+1}\,,\,\sigma_i\sigma_j=\sigma_j\sigma_i\rangle,$$
where the first group of relations ranges over $i=0,1,\cdots,n-1$ and the second group ranges over $|i-j|\geq 2$. Therefore, by (iii) of the above proposition, together with the obvious relation $\mathcal{L}_i\mathcal{L}_j=\mathcal{L}_j\mathcal{L}_i$ for $|i-j|\geq 2$, we can define an action of $A_{n+1}$ on the set of exceptional collections by $\sigma_i\cdot \mathcal{E}=\mathcal{L}_i\mathcal{E}$, for $i=0,1,\cdots,n-1$.\par\indent
For an exceptional collection $\{E_0,\cdots,E_n\}$, let $\langle E_0,\cdots,E_n\rangle$ denote the full extension-closed subcategory generated by $E_0,\cdots,E_n$.
\begin{lemma}
Let $\{E_0,\cdots,E_n\}$ be a complete Ext-exceptional collection in $\mathcal{D}$. Then $\langle E_0,\cdots,E_n\rangle$ is the heart of a bounded $t$-structure.
\end{lemma}
\begin{cor}
Let $\{E_0,\cdots,E_n\}$ be a complete Ext-exceptional collection in $\mathcal{D}$ and $(Z,\mathcal{P})$ a stability condition on $\mathcal{D}$. If $E_0,\cdots,E_n\in \mathcal{P}((0,1])$, then $\langle E_0,\cdots, E_n\rangle=\mathcal{P}((0,1])$ and $E_i$ is stable for all $i=0,1,\cdots,n$.
\end{cor}
Let $\mathcal{Q}$ denote the heart $\langle E_0,\cdots,E_n\rangle$ in the above lemma. Since $\mathcal{Q}$ is generated by $E_0,\cdots,E_n$, the Grothendieck group $K(\mathcal{Q})$ is isomorphic to the free abelian group $\mathbb{Z}^{n+1}$. In particular, a choice of complex numbers $z_0,\cdots,z_n\in H$ determines a stability function on $\mathcal{Q}$ with the Harder-Narasimhan property sending $E_i$ to $z_i$. By Proposition 2.1, this uniquely determines a locally finite stability condition on $\mathcal{D}$. \par\indent
More generally, for a complete exceptional collection $\mathcal{E}=\{E_0,\cdots,E_n\}$, we can find a sequence of integers $p=(p_0,\cdots,p_n)$ such that $\{E_0[p_0],\cdots,E_n[p_n]\}$ is a complete Ext-exceptional collection. Let $\mathcal{Q}_p=\langle E_0[p_0],\cdots, E_n[p_n]\rangle$. By the above process,
we can construct a stability condition with $\mathcal{Q}_p$ as heart by choosing $z_0,\cdots,z_n\in H$ and letting $Z(E_i[p_i])=z_i$. If the image of the central charge is a line, then we call this stability condition \emph{degenerate}; otherwise we call it \emph{nondegenerate}. Define $\Theta_{\mathcal{E}}$ as the subset of $\mathrm{Stab}(\mathcal{D})$ obtained in such way, up to the action of $\widetilde{\mathrm{GL}^+(2,\mathbb{R})}$. By Corollary 3.1, each $E_i$ is stable for any $\sigma\in \Theta_{\mathcal{E}}$. However, a stability condition for which each $E_i$ is stable need not lie in $\Theta_{\mathcal{E}}$.
\begin{lemma}
The subspace $\Theta_{\mathcal{E}}$ of $\mathrm{Stab}(\mathcal{D})$ is an open, connected and simply connected $(n+1)$-dimensional submanifold. In fact, it is homeomorphic to the space
$$C_{\mathcal{E}}=\big\{(m_0,\cdots,m_n,\phi_0,\cdots,\phi_n)\in\mathbb{R}^{2(n+1)}\,|\,m_i>0\;\textrm{for all}\;i\;\;\textrm{and}\;\;\phi_i<\phi_j+\alpha_{i,j}\;\textrm{for}\; i<j\big\},$$
where
$$\alpha_{i,j}=\min_{i<l_1<\cdots<l_s<j}\{k_{i,l_1}+k_{l_1,l_2}+\cdots+k_{l_s,j}-s\}$$
and $k_{i,j}=\min_k\{\mathrm{Hom}^k(E_i,E_j)\neq 0\}$ (if no such $k$ exists, set $k_{i,j}=+\infty$); see \cite{Mac07}, \cite{Shi13}.
\end{lemma}
The homeomorphism is given explicitly by
$$(Z,\mathcal{P})\mapsto(|Z(E_0)|,\cdots,|Z(E_n)|, \phi_{\sigma}(E_0),\cdots,\phi_{\sigma}(E_n)).$$
 The intuition for this map to be a homeomorphism is as follows. For simplicity, assume $k_{i,j}=0$ for all $i<j$. It is easy to see that $p=(p_0,\cdots,p_n)$ is a sequence of integers such that $\langle E_0[p_0],\cdots,E_n[p_n]\rangle$ is Ext-exceptional if and only if $p_j-p_i<0$ for all $i<j$. Since each $p_i$ is an integer, this is equivalent to requiring that $p_j-p_i\leq -(j-i)$ for all $i<j$. \par\indent
Fix $\sigma=(Z,\mathcal{P})\in \Theta_{\mathcal{E}}$ and $p=(p_0,\cdots,p_n)$ for which $\langle E_0[p_0],\cdots,E_n[p_n]\rangle=\mathcal{P}((0,1])$, up to the action of $\widetilde{\mathrm{GL}^+(2,\mathbb{R})}$. Thus, $\phi_{\sigma}(E_i[p_i])<\phi_{\sigma}(E_j[p_j])+1$ for all $i,j$. Therefore, we must have $\phi_{\sigma}(E_i)<\phi_{\sigma}(E_j)+(p_j-p_i+1)\leq \phi_{\sigma}(E_j)-(j-i-1)$ for all $i<j$. It is also straightforward to check the converse.\par\indent
This shows that the map is a bijection. For a proof that it is a homeomorphism, see \cite[Lemma 3.19]{Mac07}.

\section{Stability Conditions Generated by a Strong Complete Exceptional Collection}
In this section, we build up the proof of Theorem 1.1. Throughout this section, we let $\dag$ denote the following condition on an exceptional sequence $\mathcal{E}=\{E_0,\cdots,E_n\}$:
$$\textrm{(\dag)}\qquad \textrm{$\mathcal{E}$ is a strong complete exceptional sequence with no orthogonal pairs}$$
$$\textrm{such that its iterated mutations are again strong complete exceptional}.$$
This notion is similar to the notion of `geometric' or `simple' collection in \cite{Bri05}. To recall some notations, we let $S_{\mathcal{E}}$ denote that set of all iterated mutations of $\mathcal{E}$, and set $\Sigma_{\mathcal{E}}=\bigcup_{\mathcal{F}\in S_{\mathcal{E}}}\Theta_{\mathcal{E}}$. Note that by Proposition 3.1, complete exceptionality is preserved by mutations; in general, however, strongness is not preserved by mutations. \par\indent
We've already seen that the subspace $\Theta_{\mathcal{E}}$ is an open, connected and simply connected $(n+1)$-dimensional manifold. In fact, Macr\`{i} further showed that $\Sigma_{\mathcal{E}}$ is an open and connected $(n+1)$-dimensional manifold \cite[Corollary 3.20]{Mac07}. However, in order to prove that $\Sigma_{\mathcal{E}}$ is simply connected, we need to examine more closely how the $\Theta_{\mathcal{F}}$'s are glued together.
\begin{lemma}
Let $\mathcal{E}$ be $\dag$, and let $\mathcal{F}$ be a single mutation of $\mathcal{E}$. Then, $\Theta_{\mathcal{E}}\bigcap\Theta_{\mathcal{F}}$ is nonempty, path connected and simply connected.
\end{lemma}
\noindent\emph{Proof.} For nonemptiness, see \cite[Corollary 3.20]{Mac07}. \par\indent
We assume $\mathcal{F}$ is obtained from $\mathcal{E}$ by a single right mutation, i.e. $\mathcal{F}=\mathcal{R}_k\mathcal{E}=\{E_0,\cdots,E_{k-1}, E_{k+1}, \mathcal{R}_{E_{k+1}}E_k, E_{k+2}, \cdots, E_n\}$ for some $0\leq k\leq n-1$. The case of a left mutation is similar. By Lemma 3.2, since $\mathcal{E}$ is strong complete exceptional, we have $\Theta_{\mathcal{E}}\cong C_{\mathcal{E}}=$
$$\big\{(m_0,\cdots,m_n,\phi_0,\cdots,\phi_n)\in\mathbb{R}^{2(n+1)}\,|\,m_i>0\;\textrm{for all}\;i\;\;\textrm{and}\;\;\phi_i<\phi_j-(j-i-1)\;\textrm{for}\; i<j\big\},$$
where $m_i=|Z(E_i)|$ and $\phi_i=\phi(E_i)$. Therefore, $\sigma=(Z,\mathcal{P})$ lies in $\Theta_{\mathcal{E}}\bigcap\Theta_{\mathcal{F}}$ if and only if the following conditions are satisfied:\\
(i) $\phi(E_i)<\phi(E_j)-(j-i-1)$ for $i<j$; \\
(ii) $\mathcal{R}_{E_{k+1}}E_k$ is stable and $\phi(E_{k+1})<\phi(\mathcal{R}_{E_{k+1}}E_k)$;\\
(iii) $\phi(E_{k+1})<\phi(E_{k+i})-(i-1)$ for $i\geq 2$;\\
(iv) $\phi(E_{k-i})<\phi(\mathcal{R}_{E_{k+1}}E_k)-i$ for $i\geq 1$;\\
(v) $\phi(\mathcal{R}_{E_{k+1}}E_k)<\phi(E_{k+i})-(i-2)$ for $i\geq 2$. \\
We wish to express these conditions fully in terms of $\phi_i=\phi(E_i)$, for $i=0,1,\cdots,n$. To do this, we proceed in several steps.\par\indent
\emph{Step 1.} Given (i), we can replace (ii) by (ii'): $\phi(E_{k+1})<\phi(E_k)+1$. Indeed, by definition of a right mutation, we have a distinguished triangle
$$\mathrm{Hom}(E_k,E_{k+1})\otimes E_{k+1}\rightarrow\mathcal{R}_{E_{k+1}}E_k\rightarrow E_k[1].$$
Given (i), we have $E_k, E_{k+1}$ are stable by \cite[Proposition 3.17]{Mac07}. Assume $\mathcal{R}_{E_{k+1}}E_k$ is stable as well, then we must have $\phi(E_{k+1})<\phi(\mathcal{R}_{E_{k+1}}E_k)<\phi(E_k)+1$ (the maps are clearly zero, otherwise we will reach a contradiction by splitting the triangle). This shows (ii)$\Rightarrow$ (ii'). Conversely if we are given (ii'), then $\mathcal{R}_{E_{k+1}}E_k$ is stable by \cite[Proposition 3.17]{Mac07}. Then again from the above triangle we see that $\phi(E_{k+1})<\phi(\mathcal{R}_{E_{k+1}}E_k)<\phi(E_k)+1$. \par\indent
\emph{Step 2.} We show that (i), (ii),(iii) implies (iv),(v). For (iv), we note that by (i) and (ii), we have $\phi(E_{k-i})<\phi(E_{k+1})-i<\phi(\mathcal{R}_{E_{k+1}}E_k)-i$. For (v), note that in Step 1, we showed that (ii) implies $\phi(E_{k+1})<\phi(\mathcal{R}_{E_{k+1}}E_k)<\phi(E_k)+1$. Therefore, combining with (i), we have $\phi(\mathcal{R}_{E_{k+1}}E_k)<\phi(E_k)+1<\phi(E_{k+i})-(i-1)+1=\phi(E_{k+i})-(i-2)$. \par\indent
\emph{Step 3.} Finally, by Step 1 and 2, we obtain that $\Theta_{\mathcal{E}}\bigcap\Theta_{\mathcal{F}}$ is homeomorphic to $(\mathbb{R}_{>0})^{n+1}\times \Phi_{\mathcal{E}\bigcap\mathcal{F}}$, where
$$\Phi_{\mathcal{E}\bigcap\mathcal{F}}=\{(\phi_0,\cdots,\phi_n)\in\mathbb{R}^{n+1}\,|\,(i),(ii'),(iii)\}$$
by identifying $\phi_i$ with $\phi(E_i)$. It is clear that (i),(ii'),(iii) define a path connected and simply connected open submanifold of $\mathbb{R}^{n+1}$. In fact, we can easily check that this subspace is convex. Hence, we conclude that $\Theta_{\mathcal{E}}\bigcap\Theta_{\mathcal{F}}$ is path connected and simply connected. $\Box$

\begin{cor}
Let $\mathcal{E}$ be $\dag$, and let $\mathcal{F}$ be a single mutation of $\mathcal{E}$. Then, $\Theta_{\mathcal{E}}\bigcup\Theta_{\mathcal{F}}$ is path connected and simply connected.
\end{cor}
\noindent\emph{Proof.} $\Theta_{\mathcal{E}}\bigcup\Theta_{\mathcal{F}}$ is path connected since $\Theta_{\mathcal{E}}$ and $\Theta_{\mathcal{F}}$ are path connected and $\Theta_{\mathcal{E}}\bigcap\Theta_{\mathcal{F}}\neq \emptyset$. Moreover, since $\Theta_{\mathcal{E}}$ and $\Theta_{\mathcal{F}}$ are simply connected, and $\Theta_{\mathcal{E}}\bigcap\Theta_{\mathcal{F}}$ path connected, we conclude that $\Theta_{\mathcal{E}}\bigcup\Theta_{\mathcal{F}}$ is simply connected by Seifert-van Kampen theorem. $\Box$\\\par\indent
Next, we study some boundary conditions of the open subsets $\Theta_{\mathcal{E}}$. The following few lemmas follow the same idea as Lemma 4.7 through Lemma 4.11 in \cite{Mac07}.
\begin{lemma}
The closure of $\Theta_{\mathcal{E}}$ is contained in $\Sigma_{\mathcal{E}}$.
\end{lemma}
\noindent\emph{Proof.} Let $\overline{\sigma}=(\overline{Z},\overline{\mathcal{P}})$ be a stability condition in $\partial\Theta_{\mathcal{E}}$. Then, we can find integers $p_0>p_1>\cdots>p_n$ such that $\{E_0[p_0],\cdots,E_n[p_n]\}$ is Ext-exceptional and contained in $\overline{\mathcal{P}}([0,1])$. Moreover, we already know that each $E_i$ is stable in any stability condition in $\Theta_{\mathcal{E}}$, and hence they are semistable in $\overline{\sigma}$ as semistability is a closed condition. Define $N_{\overline{\sigma}}(\mathcal{E},0):=\#\{i\in \{0,1,\cdots,n\}: \phi_{\overline{\sigma}}(E_i[p_i])=0\}$ and $N_{\overline{\sigma}}(\mathcal{E},1):=\#\{i\in \{0,1,\cdots,n\}: \phi_{\overline{\sigma}}(E_i[p_i])=1\}$. Since $\overline{\sigma}\in \partial\Theta_{\mathcal{E}}$, we must have $N_{\overline{\sigma}}(\mathcal{E},0), N_{\overline{\sigma}}(\mathcal{E},1)\geq 1$. \par\indent
We first deal with the base case when $N_{\overline{\sigma}}(\mathcal{E},0)=N_{\overline{\sigma}}(\mathcal{E},1)=1$. Then, there exists a unique pair $(i,j)$ such that $\phi_{\overline{\sigma}}(E_i[p_i])=\phi_{\overline{\sigma}}(E_j[p_j])-1=0$. The exceptional collection $\mathcal{E}'=\{E_0[p_0],\cdots,E_i[p_i+1],\cdots,E_n[p_n]\}$ is contained in $\overline{\mathcal{P}}((0,1])$, and hence is not Ext (otherwise $\overline{\sigma}\in\Theta_{\mathcal{E}}$). Since $p_0>p_1>\cdots>p_n$ are integers, this implies that $p_i+1=p_{i-1}$ and thus $\mathrm{Hom}(E_{i-1}[p_{i-1}],E_i[p_i+1])\cong \mathrm{Hom}(E_{i-1},E_i)\neq 0$. We discuss two cases. \par\indent
1) If $j\neq i-1$, then in particular, $E_{i-1}[p_{i-1}]\in \overline{\mathcal{P}}((0,1))$. By the exact triangle
$$E_{i-1}[p_i]\rightarrow \mathrm{Hom}(E_{i-1},E_i)\otimes E_i[p_i]\rightarrow \mathcal{R}_{E_i}E_{i-1}[p_i]\rightarrow E_{i-1}[p_{i-1}],$$
we deduce that $\mathcal{R}_{E_i}E_{i-1}[p_i]\in \overline{\mathcal{P}}((0,1))$. Therefore, the collection
$$\{E_0[p_0],\cdots, E_{i-2}[p_{i-2}], E_i[p_{i-1}], \mathcal{R}_{E_i}E_{i-1}[p_i], E_{i+1}[p_{i+1}], \cdots, E_n[p_n]\}$$
is Ext-exceptional and contained in $\overline{\mathcal{P}}((0,1])$, which implies that $\overline{\sigma}\in\Theta_{\mathcal{R}_{i-1}\mathcal{E}}$. \par\indent
2) If $j=i-1$, we do the following. Let $\sigma_s\rightarrow \overline{\sigma}$ such that each $\sigma_s\in\Theta_{\mathcal{E}}$ has $\langle E_0[p_0],\cdots,E_n[p_n]\rangle$ as heart. Since $\phi_{\overline{\sigma}}(E_i[p_i])=0$, we have $\phi_s(E_i[p_i])<\phi_s(E_l[p_l])$ for all $l\neq i$ and $s>>0$. By \cite[Proposition 3.17]{Mac07}, $\sigma_s$ induces a stability condition on $\mathrm{Tr}(E_{i-1}, E_i)$. Let $m=\dim \mathrm{Hom}(E_{i-1}, E_i)$, then we know that $\mathrm{Tr}(E_{i-1}, E_i)\cong D^b(P_m)$, where $P_m$ is the quiver with two vertices and $m$ arrows. Hence, $\overline{\sigma}$ induces a stability condition on $D^b(P_m)$. By \cite[Lemma 4.2]{Mac07}, there exists a stable exceptional pair $(F_{i-1}, F_i)$ consisting of shifts of an iterated mutation of $(E_{i-1}, E_i)$, such that $\{E_0[p_0],\cdots, E_{i-2}[p_{i-2}], F_{i-1}, F_i, E_{i+1}[p_{i+1}],\cdots, E_n[p_n]\}\subset \overline{\mathcal{P}}((0,1])$ is Ext. \par\indent
For general $\overline{\sigma}$, we first find a sequence $\sigma_s\rightarrow \overline{\sigma}$ with $\sigma_s\in \overline{\Theta}_{\mathcal{E}}$ satisfying $N_{\sigma_s}(\mathcal{E},0)=N_{\sigma_s}(\mathcal{E},1)=1$ for all $s$. By the previous paragraph, we may find some $l_1\in \{\mathcal{L}_0,\cdots,\mathcal{L}_{n-1},\mathcal{R}_0,\cdots,\mathcal{R}_{n-1}\}$ and integer $k_1$ such that $\sigma_s\in \Theta_{l_1^{k_1}\mathcal{E}}$ for all $s$. If $\overline{\sigma}\in \Theta_{l_1^{k_1}\mathcal{E}}$, we are done. Assume otherwise that $\overline{\sigma}\in \partial\Theta_{l_1^{k_1}\mathcal{E}}$. Then, we can repeat the above argument by deforming $\overline{\sigma}$ in $\overline{\Theta}_{l_1^{k_1}\mathcal{E}}$, and obtain a new $l_2$ such that $\overline{\sigma}\in\overline{\Theta}_{l_2^{k_2}l_1^{k_1}\mathcal{E}}$ for some integer $k_2$. I claim that this process will terminate. Indeed, at each step we are constructing a Jordan-Holder filtration of the $E_i$'s by semistable objects of the same phase. However, as $\overline{\sigma}$ is locally finite, this process must end within finite steps.$\Box$
\begin{cor}
Let $\mathcal{E}$ be $\dag$ and $\mathcal{F}$ an iterated mutation of $\mathcal{E}$. Assume that there exists $\overline{\sigma}=(\overline{Z},\overline{P})\in \partial\Theta_{\mathcal{E}}\bigcap\Theta_{\mathcal{F}}$ such that the image of $\overline{Z}$ is contained in a line. Then there exists $l=l_s\cdots l_1$ with $l_i\in \{\mathcal{L}_0,\cdots,\mathcal{L}_{n-1},\mathcal{R}_0,\cdots,\mathcal{R}_{n-1}\}$ for all $i$, such that $\mathcal{F}=l\mathcal{E}$, and real numbers
$0=a_0<a_1<\cdots<a_s<a_{s+1}=1$ and a continuous path $\gamma: [0,1]\rightarrow \Sigma_{\mathcal{E}}$ such that $\gamma([a_k,a_{k+1}))\subset \Theta_{l_k\cdots l_1\mathcal{E}}$ for all $k=0,1,\cdots,s$ and $\gamma(1)=\overline{\sigma}$.
\end{cor}
\noindent\emph{Proof.} By Lemma 4.2, there exists a sequence $l_1,\cdots, l_s$ (with potential repetitions) with $l_i\in \{\mathcal{L}_0,\cdots,\mathcal{L}_{n-1},\mathcal{R}_0,\cdots,\mathcal{R}_{n-1}\}$ for all $i$ and $\overline{\sigma}\in \partial\Theta_{\mathcal{E}}\bigcap\Theta_{l_s\cdots l_1\mathcal{E}}$. I claim that $l_s\cdots l_1\mathcal{E}=\mathcal{F}$. Indeed, since $\overline{\sigma}\in \Theta_{\mathcal{F}}$ is degenerate, the objects of $\mathcal{F}$ are the only stable objects in $\overline{\mathcal{P}}((0,1])$. The same holds for $l_s\cdots l_1\mathcal{E}$, and thus it must agree with $\mathcal{F}$. The statement then follows from Lemma 4.1. $\Box$\\\\
The following lemma states that any loop in $\Sigma_{\mathcal{E}}$ can be decomposed into a sequence of segments such that two adjacent segments `differ' by a single mutation.
\begin{prop}
Let $\gamma: [0,1]\rightarrow \Sigma_{\mathcal{E}}$ be a continuous loop with $\gamma(0)=\gamma(1)\in \Theta_{\mathcal{E}}$. Then, up to replacing $\gamma$ by a homotopic path, there exists $l=l_s\cdots l_1$ with $l_i\in \{\mathcal{L}_0,\cdots,\mathcal{L}_{n-1},\mathcal{R}_0,\cdots,\mathcal{R}_{n-1}\}$ for all $i$, such that $l\mathcal{E}=\mathcal{E}$, and real numbers $0=a_0<a_1<\cdots<a_s<a_{s+1}=1$ such that $\gamma([a_k,a_{k+1}))\subset \Theta_{l_k\cdots l_1\mathcal{E}}$ for all $k=0,1,\cdots,s$.
\end{prop}
Before proving this proposition, we first prove two lemmas.
\begin{lemma}
Let $\mathcal{F}$ be an iterated mutation of $\mathcal{E}$ and let $\gamma: [0,1]\rightarrow \Theta_{\mathcal{E}}\bigcup\Theta_{\mathcal{F}}$ be a continuous path such that $\gamma([0,1))\in \Theta_{\mathcal{E}}$ and $\gamma(1)\in \partial\Theta_{\mathcal{E}}\bigcap \Theta_{\mathcal{F}}$. Then there exist $\gamma'$ with $\gamma'([0,1))\in \Theta_{\mathcal{E}}$ and $\gamma'(1)\in \partial\Theta_{\mathcal{E}}\bigcap \Theta_{\mathcal{F}}$ degenerate for $\mathcal{F}$, and some $\gamma''\subset \Theta_{\mathcal{F}}$, such that $\gamma$ is homotopic to $\gamma''\circ \gamma'$.
\end{lemma}
\noindent\emph{Proof.} There exists a continuous sequence $G_s\in \widetilde{\mathrm{GL}^+(2,\mathbb{R})}$ (setting $G_0=\mathrm{id}$) such that $G_s\cdot \gamma(1)\rightarrow \overline{\sigma}$, where $\overline{\sigma}\in \Theta_{\mathcal{F}}$ is a degenerate stability condition. It is clear that $\overline{\sigma}\in \partial\Theta_{\mathcal{E}}$. Let $\gamma''$ denote the path $G_s\cdot \gamma(1)\rightarrow \overline{\sigma}$, which is contained in $\overline{\Theta}_{\mathcal{E}}\bigcap\Theta_{\mathcal{F}}$. Find any path $\gamma'$ such that $\gamma'(0)=\gamma(0), \gamma'([0,1))\subset \Theta_{\mathcal{E}}$ and $\gamma'(1)=\overline{\sigma}$, which is possible since $\Theta_{\mathcal{E}}$ is path connected and $\Sigma_{\mathcal{E}}$ is locally Euclidean. Since $\Theta_{\mathcal{E}}$ is simply connected and $\Sigma_{\mathcal{E}}$ is locally Euclidean, we conclude that $\gamma\sim \gamma''\circ\gamma'$. $\Box$
\begin{lemma}
Let $\mathcal{F}$ be an iterated mutation of $\mathcal{E}$ and let $\gamma: [0,1]\rightarrow \Theta_{\mathcal{E}}\bigcup\Theta_{\mathcal{F}}$ be a continuous path such that $\gamma([0,1))\in \Theta_{\mathcal{E}}$ and $\gamma(1)\in \partial\Theta_{\mathcal{E}}\bigcap \Theta_{\mathcal{F}}$ is degenerate. Then, up to replacing $\gamma$ by a homotopic path, there exists $l=l_s\cdots l_1$ with $l_i\in \{\mathcal{L}_0,\cdots,\mathcal{L}_{n-1},\mathcal{R}_0,\cdots,\mathcal{R}_{n-1}\}$ for all $i$, such that $\mathcal{F}=l\mathcal{E}$, and real numbers
$0=a_0<a_1<\cdots<a_s<a_{s+1}=1$ such that $\gamma([a_k,a_{k+1}))\subset \Theta_{l_k\cdots l_1\mathcal{E}}$ for all $k=0,1,\cdots,s$.
\end{lemma}
\begin{figure}[H]
 \centering
 \includegraphics[height=11.5cm, width=14cm]{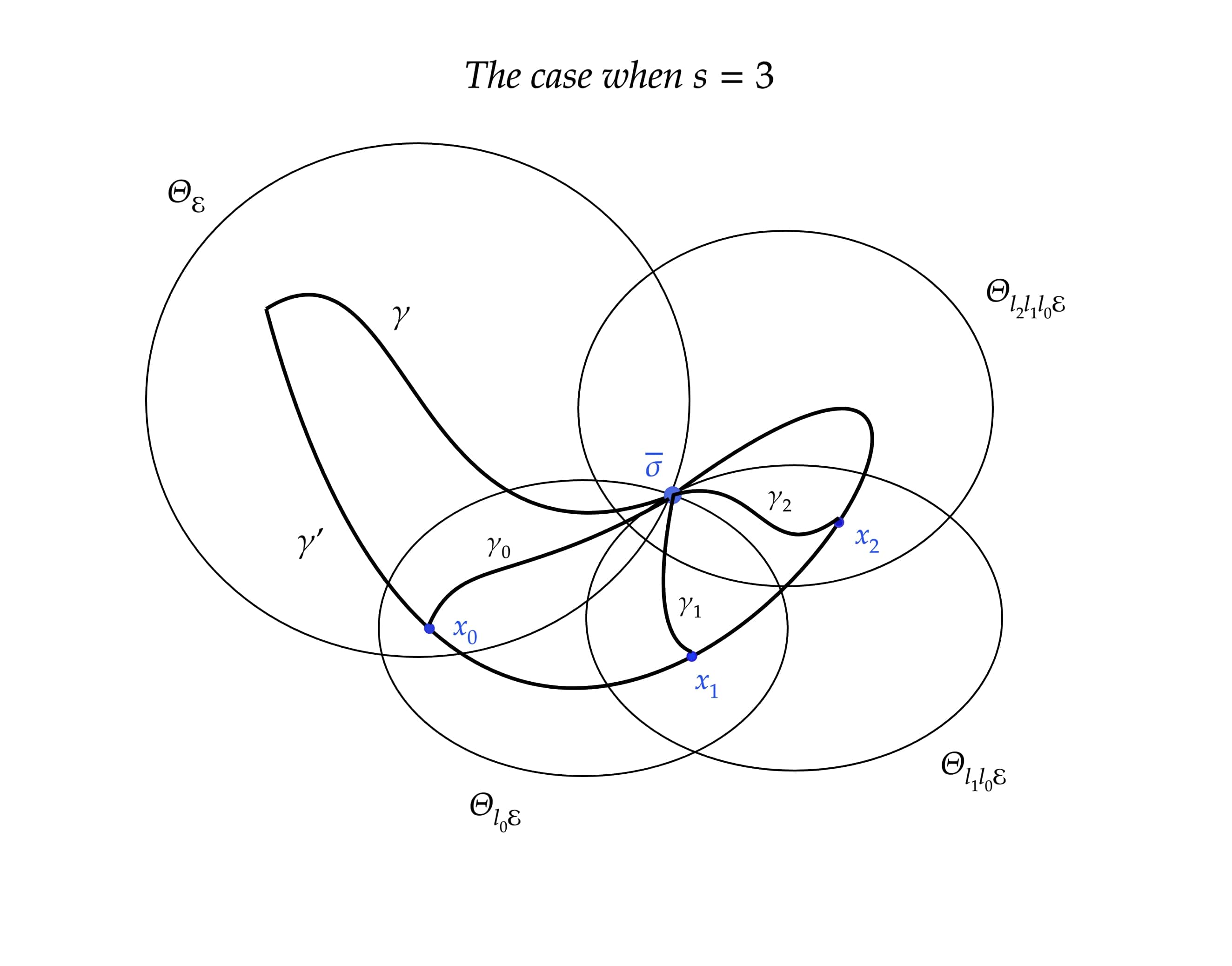}
\end{figure}
\noindent\emph{Proof.} By Corollary 4.2, we can find $l=l_s\cdots l_1$ with $\mathcal{F}=l\mathcal{E}$, real numbers $0=a_0<a_1<\cdots<a_s<a_{s+1}=1$ and a continuous path $\gamma': [0,1]\rightarrow \Sigma_{\mathcal{E}}$ such that $\gamma'([a_k,a_{k+1}))\subset \Theta_{l_k\cdots l_1\mathcal{E}}$ for all $k=0,1,\cdots,s$ and $\gamma'(1)=\overline{\sigma}$. Since $\Theta_{\mathcal{E}}$ is path connected, we may assume that $\gamma'(0)=\gamma(0)$. We need to show that $\gamma$ and $\gamma'$ are homotopic. \par\indent
Notice that in the proof of Corollary 4.2, the sequence $l_1,\cdots, l_s$ are chosen such that $\overline{\sigma}\in \overline{\Theta}_{l_i\cdots l_1\mathcal{E}}$ for all $i$, and we terminate at the first $s$ such that $\overline{\sigma}\in \Theta_{l_s\cdots l_1\mathcal{E}}$. In particular, $\overline{\sigma}\in (\bigcap_{i=0}^{s-1} \partial \Theta_{l_i\cdots l_1})\bigcap \Theta_{l_s\cdots l_1\mathcal{E}}$. By Lemma 4.1 and the fact that $\Sigma_{\mathcal{E}}$ is locally Euclidean, for each $i=0,1,\cdots,s-1$, we can find a point $x_i\in \Theta_{l_i\cdots l_1\mathcal{E}}\bigcap \Theta_{l_{i+1}\cdots l_1\mathcal{E}}\bigcap \gamma'$ and a path $\gamma_i$ such that $\gamma_i(0)=x_i, \gamma_i([0,1))\subset \Theta_{l_i\cdots l_1\mathcal{E}}\bigcap \Theta_{l_{i+1}\cdots l_1\mathcal{E}}$ and $\gamma_i(1)=\overline{\sigma}$. This shows that $\gamma$ is homotopic to $\gamma'$, as illustrated in the above figure. $\Box$\\\\
\noindent\emph{Proof of Proposition 4.1.} Let $\gamma: [0,1]\rightarrow \Sigma_{\mathcal{E}}$ be a continuous loop. Without loss of generality, we may assume $\gamma(0)=\gamma(1)\in \Theta_{\mathcal{E}}$. Since $\gamma$ is compact, we may find real numbers $0=a_0<a_1<\cdots<a_m<a_{m+1}=1$ and $\mathcal{F}_0=\mathcal{E},\mathcal{F}_1,\cdots,\mathcal{F}_{m-1}, \mathcal{F}_m=\mathcal{E}\in S_{\mathcal{E}}$ such that $\gamma([a_k,a_{k+1}))\subset \Theta_{\mathcal{F}_k}$ for all $k=0,1,\cdots,m$. By Lemma 4.3, we may assume that $\gamma(a_k)\in \partial\Theta_{\mathcal{F}_{k-1}}\bigcap\Theta_{\mathcal{F}_k}$ is degenerate in $\mathcal{F}_k$, for all $k=1,2,\cdots,m$. \par\indent
Define $\gamma_k$ as $\gamma([a_k,a_{k+1}])$ for $k=0,1,\cdots,m-1$. By Lemma 4.4, we can find $l_k=l_{s_k,k}\cdots l_{1,k}$ and real numbers $0=b_{0,k}<\cdots<b_{s_k,k}<b_{s_k+1,k}=1$ such that $\mathcal{F}_{k+1}=l_k\mathcal{F}_k$ and $\gamma_k([b_{i,k},b_{i+1,k}))\subset \Theta_{l_{i,k}\cdots l_{1,k}\mathcal{F}_k}$ for all $i=0,1,\cdots,s_k$. By combining these as $k$ ranges over $1,2,\cdots,m$, we conclude the proof of Proposition 4.1. $\Box$
\begin{cor}
With notation as above. If $\bigcap_{\mathcal{F}\in S_{\mathcal{E}}}\Theta_{\mathcal{E}}\neq \emptyset$, then $\Sigma_{\mathcal{E}}$ is simply connected.
\end{cor}
In particular, this implies that $\mathrm{Stab}(\mathbb{P}^1)$ is simply connected.
\section{Simply Connectedness in the Case of $\mathbb{P}^3$}
In this section, we apply results in the previous section to show the following theorem.
\begin{thm}
Let $\mathcal{E}$ be the exceptional collection $\{\mathcal{O},\mathcal{O}(1),\mathcal{O}(2),\mathcal{O}(3)\}$ on $\mathbb{P}^3$. Then, the subspace $\Sigma_{\mathcal{E}}$ of $\mathrm{Stab}(\mathbb{P}^3)$ is simply connected.
\end{thm}
The idea of the proof is as follows. \par\indent
Using Proposition 4.1, we can associate to any continuous loop $\gamma$ in $\Sigma_{\mathcal{E}}$ a `pattern', which is a word $l=l_s\cdots l_1$ with $l_i\in \{\mathcal{L}_0,\cdots,\mathcal{L}_{n-1},\mathcal{R}_0,\cdots,\mathcal{R}_{n-1}\}$ for all $i$. Suppose that the action of $A_4$ on $S_{\mathcal{E}}$ is free. This implies that $l=1$ and thus $l$ must be a combination of the relations of $A_4$. Hence, it suffices to check that any loop whose pattern is one of the relations is contractible, which can be done by some straightforward calculations. \par\indent
Now we prove that $A_4$ acts freely on the set of iterated mutations of $\mathcal{E}$. First, let's recall the following basic facts about the braid group $A_4$.
\begin{lemma}
1) The center of $A_4$ is generated by $(\sigma_0\sigma_1\sigma_2)^4=(\sigma_2\sigma_1\sigma_0)^4$.\\
2) The element $\delta=(\sigma_0\sigma_1\sigma_2)(\sigma_0\sigma_1)\sigma_0\in A_4$ has the property that $\delta^{-1}\sigma_i\delta=\sigma_{2-i}$ for $i=0,1,2$.
\end{lemma}
\noindent\emph{Proof.} See \cite[Lemma 2.1]{Bri05}. $\Box$\\\par\indent
Let $\mathcal{E}=\{E_0,E_1,E_2,E_3\}$ be an exceptional collection on $\mathbb{P}^3$. We define $\mathcal{F}=\delta \mathcal{E}$ to be the \emph{left dual collection} of $\mathcal{E}$. As an example, the left dual collection of $\{\mathcal{O},\mathcal{O}(1),\mathcal{O}(2),\mathcal{O}(3)\}$ is $\{\Omega^3(3),\Omega^2(2),\Omega^1(1),\mathcal{O}\}$, where $\Omega^k$ is the sheaf of holomorphic $k$-forms on $\mathbb{P}^3$. For some technical reasons we will see later, it is often more convenient to consider the braid action after passing to the dual collection. Finally, as an abuse of notation, we will use $\mathcal{L}_i\,(\mathcal{R}_i)$ interchangeably with $\sigma_i\,(\sigma_i^{-1})$ for the rest of the section.\par\indent
Recall that the Euler form on Grothendieck group $K(\mathbb{P}^3)$ is given by $\chi([E],[F])=\sum_{i}(-1)^i\dim\mathrm{Hom}^i(E,F)$. Under the basis $\{[E_i]\}$ of $K(\mathbb{P}^3)$ given by the exceptional collection, the Gram matrix of $\chi$ is given by $A=(a_{ij})$, where $a_{ij}=\dim\mathrm{Hom}(E_i,E_j)$ for $i,j=0,1,2,3$. By the strong exceptionality of $\mathcal{E}$, $A$ in fact takes the form of an upper triangular matrix
\[
\begin{pmatrix}
1&a_{01}&a_{02}&a_{03}\\
0&1&a_{12}&a_{13}\\
0&0&1&a_{23}\\
0&0&0&1
\end{pmatrix}.
\]
By Serre duality, we have $\chi([E],[F])=\chi([F],[\kappa E])$, where $\kappa=-\otimes \mathcal{O}(-4)[3]$ is the Serre functor on $D^b(\mathbb{P}^3)$. From this we deduce that $\kappa=A^{-1}A^T$ under the above basis. By [BP94, Lemma 3.1], $-\kappa$ is unipotent. Explicit computation via $\kappa=A^{-1}A^T$ shows that the unipotency of $-\kappa$ is equivalent to the following two conditions:
\begin{equation}
a_{01}^2+a_{02}^2+a_{03}^2+a_{12}^2+a_{13}^2+a_{23}^2-a_{01}a_{12}a_{02}-a_{01}a_{13}a_{03}-a_{02}a_{23}a_{03}-a_{12}a_{23}a_{13}+a_{01}a_{12}a_{23}a_{03}-8=0.
\end{equation}
\begin{equation}
a_{01}^2a_{23}^2+a_{02}^2a_{12}^2+a_{03}^2a_{12}^2-2a_{01}a_{13}a_{02}a_{23}-2a_{02}a_{12}a_{03}a_{13}+2a_{01}a_{12}a_{23}a_{03}-16=0.
\end{equation}
Let $\Gamma\subset \mathbb{Z}^6$ denote the set of integer six-tuples satisfying (1) and (2). Then, there is a map $T: S_{\mathcal{E}}\rightarrow \Gamma$ defined by $$\{F_0,F_1,F_2,F_3\}\mapsto (a_{01},a_{02},a_{03},a_{12},a_{13},a_{23}),$$
where $a_{ij}=\dim\mathrm{Hom}(F_i,F_j)$. \par\indent
Now, we define a group $G$ by its generators and relations:
$$G:=\langle v, w_2,w_3\,|\,v^2=w_2^4=w_3^2=1, w_3^{-1}w_2w_3=w_2^{-1}, (vw_3w_2^3)^3=1, (vw_3vw_2^2)^2=1\rangle.$$
In particular, the subgroup of $G$ generated by $w_2,w_3$ is isomorphic to the dihedral group $D_8$. Moreover, $\Gamma$ carries a $G$-action given by
\begin{align*}
v&: (a_{01},a_{02},a_{03},a_{12},a_{13},a_{23})\mapsto(a_{01},a_{03},a_{02},a_{01}a_{03}-a_{13},a_{01}a_{02}-a_{12},a_{23}),\\
w_2&: (a_{01},a_{02},a_{03},a_{12},a_{13},a_{23})\mapsto (a_{03},a_{13},a_{23},a_{01},a_{02},a_{12}),\\
w_3&:(a_{01},a_{02},a_{03},a_{12},a_{13},a_{23})\mapsto(a_{23},a_{13},a_{03},a_{12},a_{02},a_{01}).
\end{align*}
\begin{lemma}
The map $f: A_4\rightarrow G$ given by
$$\mathcal{R}_0\mapsto w_2^2vw_3,\quad \mathcal{R}_1\mapsto w_2vw_3w_2,\quad \mathcal{R}_2\mapsto vw_3w_2^2$$
is a group homomorphism. Moreover, $T\delta$ is equivariant with respect to $f$, in the sense that
$$T\delta(\sigma(\mathcal{F}))=f(\sigma)T\delta(\mathcal{F})$$
for all $\sigma\in A_4$ and all $\mathcal{F}\in S_{\mathcal{E}}$.
\end{lemma}
\noindent\emph{Proof.} To show that $f$ is a group homomorphism, it suffices to show that it sends all relations in $A_4$ to the identity in $G$. The computations are straightforward. For example,
$$f(\mathcal{R}_0\mathcal{R}_1\mathcal{R}_0)=(w_2^2vw_3)(w_2vw_3w_2^2)(w_2vw_3)=(w_2^2vw_3w_2)(w_2w_3vw_2w_3v)(vw_3)$$
$$=w_2^2vw_2^2vw_2=w_2^2w_3vw_2^2vw_3w_2=w_2(w_2w_3vw_2w_3v)vw_3w_2vw_3w_2=f(\mathcal{R}_1\mathcal{R}_0\mathcal{R}_1).$$
The other two relations can be checked similarly. \par\indent
To show that $T\delta$ is equivariant, it suffices to check $T\delta(\sigma(\mathcal{F}))=f(\sigma)T\delta(\mathcal{F})$ for $\sigma=\mathcal{R}_0,\mathcal{R}_1,\mathcal{R}_2$. We will first assume $\sigma=\mathcal{R}_2$. The defining triangle
$$F_0\rightarrow \mathrm{Hom}(F_0,F_1)\otimes F_1\rightarrow \mathcal{R}_{F_1}F_0\rightarrow F_0[1]$$
and the strong exceptionality of $\mathcal{R}_0\mathcal{F}=\{F_1,\mathcal{R}_{F_1}F_0, F_2,F_3\}$ implies that $T(\mathcal{R}_0\mathcal{F})=f(\mathcal{R}_2)T(\mathcal{F})$. Therefore,
$$T\delta(\mathcal{R}_2\mathcal{F})=T\mathcal{R}_0(\delta\mathcal{F})=f(\mathcal{R}_2)T\delta(\mathcal{F}).$$
The proof for $\sigma=\mathcal{R}_0$ and $\sigma=\mathcal{R}_1$ are similar.$\Box$
\begin{prop}
Let $\mathcal{E}=\{\mathcal{O},\mathcal{O}(1),\mathcal{O}(2),\mathcal{O}(3)\}$. Then the action of $A_4$ on $S_{\mathcal{E}}$ is free.
\end{prop}
\noindent\emph{Proof.} The dual collection of $\mathcal{E}$ is $\delta\mathcal{E}=\{\Omega^3(3),\Omega^2(2),\Omega^1(1),\mathcal{O}\}$, with $T(\delta\mathcal{E})=(4,6,4,4,6,4)$. The stabilizer subgroup of this element is $G_{(4,6,4,4,6,4)}=\langle w_2,w_3\rangle\subset G$. Next, from $f(\mathcal{R}_2\mathcal{R}_1\mathcal{R}_0)=w_2$ we deduce that $\mathrm{im}f=\langle w_2, vw_3\rangle \subset G$.\par\indent
I claim that $\ker f=\langle (\mathcal{R}_2\mathcal{R}_1\mathcal{R}_0)^4\rangle=Z(A_4)$. To prove this, it suffices to show that $A_4/Z(A_4)$ is isomorphic to $\mathrm{im}f$, which is isomorphic to the abstract two-generator group given by $G'=\langle v', w_2: w_2^4=1,v'^2w_2^2v'^{-2}=w_2^2, (v'w_2^3)^3=1\rangle$ (under the obvious identification $vw_3\leftrightarrow v', w_2\leftrightarrow w_2$). Then, the map $f: A_4\rightarrow G$ is equivalent to $f': A_4\rightarrow G'$ defined by $\mathcal{R}_0\mapsto w_2^2v', \mathcal{R}_1\mapsto w_2v'w_2,\mathcal{R}_2\mapsto v'w_2^2$. However, $f'$ has an explicit inverse $f'^{-1}: G'\rightarrow A_4/Z(A_4)$ given by $v'\mapsto (\mathcal{R}_2\mathcal{R}_1\mathcal{R}_0)^2\mathcal{R}_0, w_2\mapsto \mathcal{R}_2\mathcal{R}_1\mathcal{R}_0$. This shows that $A_4/Z(A_4)\cong \mathrm{im} f$. \par\indent
Assume that some $\sigma\in A_4$ fixes $\mathcal{E}$. In particular, $f(\sigma)$ must also fix $T\delta(\mathcal{E})=(4,6,4,4,6,4)$. Therefore, $f(\sigma)\in\mathrm{im}f\bigcap G_{(4,6,4,4,6,4)}=\langle w_2,vw_3\rangle\bigcap \langle w_2, w_3\rangle=\langle w_2\rangle$. Since $\ker f=\langle (\mathcal{R}_2\mathcal{R}_1\mathcal{R}_0)^4\rangle$, this implies that $\sigma=(\mathcal{R}_2\mathcal{R}_1\mathcal{R}_0)^k$ for some integer $k$. \par\indent
\cite[Theorem 4.1]{Bon90} tells us that $(\mathcal{R}_2\mathcal{R}_1\mathcal{R}_0)^4$ acts on an exceptional collection by twisting by the anticanonical bundle. Thus, $\sigma^4=(\mathcal{R}_2\mathcal{R}_1\mathcal{R}_0)^{4k}=-\otimes \mathcal{O}(4k)$. By assumption, however, $\sigma^4$ also fixes $\mathcal{E}$, which is clearly impossible unless $k=0$. Hence, the stabilizer subgroup of $\mathcal{E}=\{\mathcal{O},\mathcal{O}(1),\mathcal{O}(2),\mathcal{O}(3)\}$ is trivial. Since $A_4$ acts transitively on $S_{\mathcal{E}}$, and in particular, all stabilizer subgroups are conjugate to each other, we conclude that the action of $A_4$ on $S_{\mathcal{E}}$ is free. $\Box$\\\\
\emph{Proof of Theorem 5.1.} Let $\mathcal{E}=\{\mathcal{O},\mathcal{O}(1),\mathcal{O}(2),\mathcal{O}(3)\}$, and fix a continuous loop $\gamma: [0,1]\rightarrow \Sigma_{\mathcal{E}}$. By Proposition 4.1, there exists $l=l_s\cdots l_1$ with $l_i\in \{\mathcal{L}_0,\mathcal{L}_1,\mathcal{L}_{2},\mathcal{R}_0,\mathcal{R}_1,\mathcal{R}_{2}\}$ for all $i$, such that $l\mathcal{E}=\mathcal{E}$, and real numbers $0=a_0<a_1<\cdots<a_s<a_{s+1}=1$ such that $\gamma([a_k,a_{k+1}))\subset \Theta_{l_k\cdots l_1\mathcal{E}}$ for all $k=0,1,\cdots,s$. Since the action of $A_4$ on $S_{\mathcal{E}}$ is free, we must have
$$l=(h_1 r_1^{\pm1} h_1^{-1})(h_2 r_2^{\pm 1} h_2^{-1})\cdots(h_s r_s^{\pm 1} h_s^{-1}),$$
where each $h_i$ is a word in $\{\mathcal{L}_0,\mathcal{L}_1,\mathcal{L}_{2},\mathcal{R}_0,\mathcal{R}_1,\mathcal{R}_{2}\}$, and each $r_i$ is one of $\mathcal{R}_0\mathcal{R}_1\mathcal{R}_0\mathcal{L}_1\mathcal{L}_0\mathcal{L}_1$, $\mathcal{R}_1\mathcal{R}_2\mathcal{R}_1\mathcal{L}_2\mathcal{L}_1\mathcal{L}_2$ or $ \mathcal{R}_0\mathcal{R}_2\mathcal{L}_0\mathcal{L}_2$. It is clear that any loop with pattern of the form $hh^{-1}$, where $h$ is a word in $\{\mathcal{L}_0,\mathcal{L}_1,\mathcal{L}_{2},\mathcal{R}_0,\mathcal{R}_1,\mathcal{R}_{2}\}$, is contractible. Therefore, it remains to show that any continuous loop with pattern one of $\{\mathcal{R}_0\mathcal{R}_1\mathcal{R}_0\mathcal{L}_1\mathcal{L}_0\mathcal{L}_1, \mathcal{R}_1\mathcal{R}_2\mathcal{R}_1\mathcal{L}_2\mathcal{L}_1\mathcal{L}_2, \mathcal{R}_0\mathcal{R}_2\mathcal{L}_0\mathcal{L}_2\}$ is contractible. \par\indent
If $l=\mathcal{R}_i\mathcal{R}_{i+1}\mathcal{R}_i\mathcal{L}_{i+1}\mathcal{L}_i\mathcal{L}_{i+1}, i=0,1$, the proof is the same as in the case of $\mathbb{P}^2$ proved in \cite[Lemma 7.8]{BM11}. So we consider the case $l=\mathcal{R}_0\mathcal{R}_2\mathcal{L}_0\mathcal{L}_2$. By assumption, $\gamma$ runs through the regions
$$\Theta_{\mathcal{E}}\rightarrow\Theta_{\mathcal{L}_2\mathcal{E}}\rightarrow\Theta_{\mathcal{L}_0\mathcal{L}_2\mathcal{E}}\rightarrow \Theta_{\mathcal{R}_2\mathcal{L}_0\mathcal{L}_2\mathcal{E}}\rightarrow\Theta_{\mathcal{E}}.$$
Let $\mathcal{F}=\{F_0,F_1,F_2,F_3\}$ be any collection satisfying $\dag$. We wish to show that $\Theta_{\mathcal{F}}\bigcap\Theta_{\mathcal{L}_2\mathcal{F}}\bigcap\Theta_{\mathcal{L}_0\mathcal{L}_2\mathcal{F}}\neq \emptyset$. By a similar argument as in Lemma 4.1, a stability condition $\sigma=(Z,\mathcal{P})$ lies in $\Theta_{\mathcal{F}}\bigcap\Theta_{\mathcal{L}_2\mathcal{F}}\bigcap\Theta_{\mathcal{L}_0\mathcal{L}_2\mathcal{F}}$ if and only if the following conditions hold:
\begin{center}
(i) $\phi(F_i)<\phi(F_j)-(j-i-1)$ for $i<j$;\;(ii) $\phi(F_0)<\phi(F_2)-2$ and $\phi(F_1)<\phi(F_2)-1$;\;(iii) $\phi(F_1)-1<\phi(F_0)$ and $\phi(F_3)-1<\phi(F_2)$.
\end{center}
These constraints clearly defines a nonempty subspace of $\mathrm{Stab}(\mathbb{P}^3)$. A similar argument shows that $\Theta_{\mathcal{F}}\bigcap\Theta_{\mathcal{R}_2\mathcal{F}}\bigcap\Theta_{\mathcal{R}_0\mathcal{R}_2\mathcal{F}}\neq \emptyset$. \par\indent
This implies that up to homotopy, we can assume that $\gamma$ lies in $\Theta_{\mathcal{E}}\bigcup\Theta_{\mathcal{L}_0\mathcal{L}_2\mathcal{E}}$, as illustrated in the figure below.
\begin{figure}[H]
 \centering
 \includegraphics[height=10cm, width=14cm]{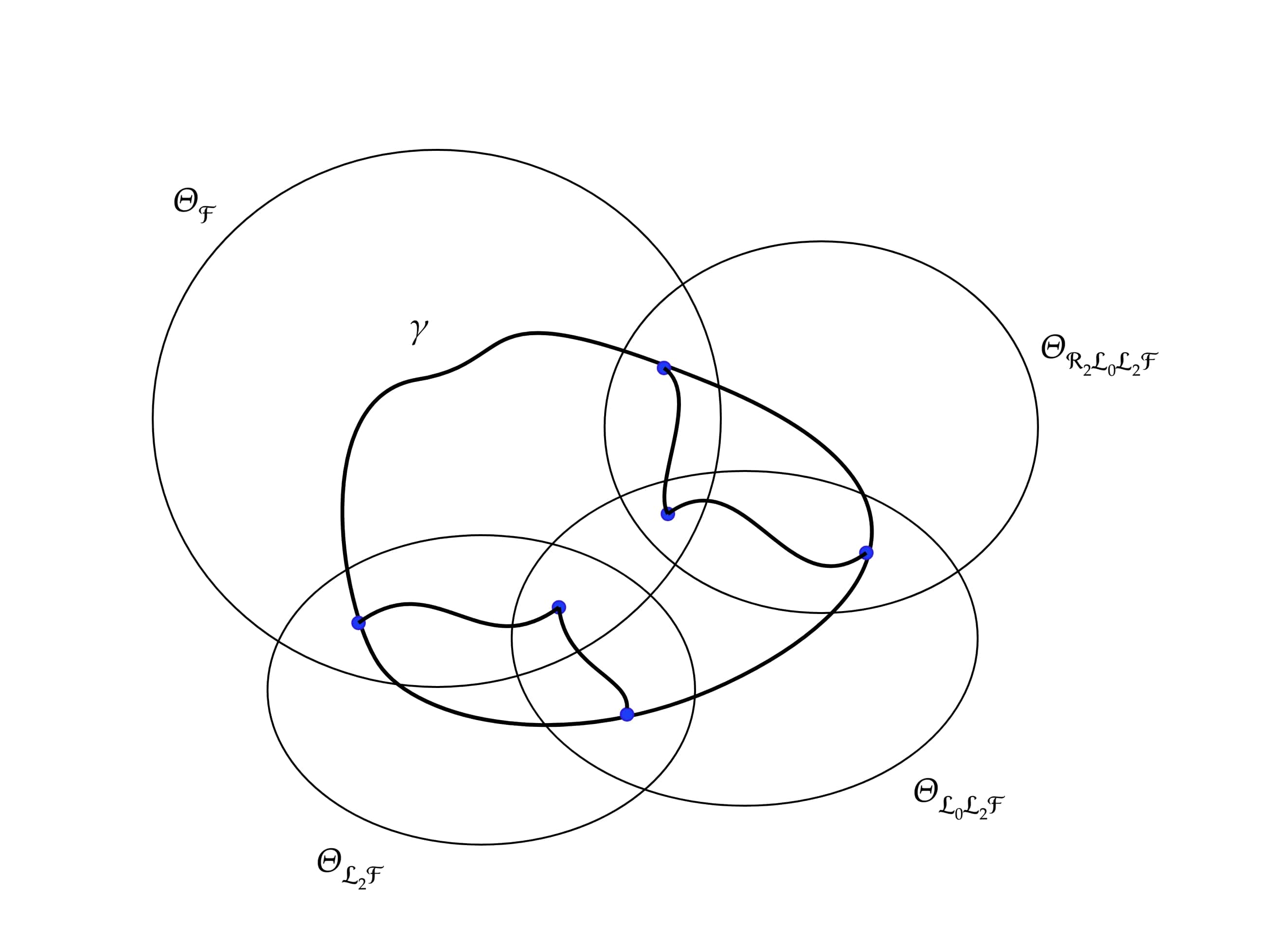}
\end{figure}
Therefore, we just need to show that $\Theta_{\mathcal{E}}\bigcup\Theta_{\mathcal{L}_0\mathcal{L}_2\mathcal{E}}$ is simply connected. By Siefert-van Kampen theorem, it is sufficient to show that $\Theta_{\mathcal{E}}\bigcap\Theta_{\mathcal{L}_0\mathcal{L}_2\mathcal{E}}$ is path connected. However, the region $\Theta_{\mathcal{E}}\bigcap\Theta_{\mathcal{L}_0\mathcal{L}_2\mathcal{E}}$ corresponds to the loci
\begin{center}
(i) $\phi(E_i)<\phi(E_j)-(j-i-1)$ for $i<j$;\;(ii) $\phi(E_1)-1<\phi(E_0)$ and $\phi(E_3)-1<\phi(E_2)$;\;(iii) $\phi(\mathcal{L}_{E_0}E_1)<\phi(E_2)-2$,
\end{center}
which is clearly path connected. This concludes the proof of the theorem. $\Box$

\section*{Appendix}
\subsection*{Some homological algebra}
We first briefly review some basics of derived functors and derived categories. For a more detailed introduction, we refer the readers to Chapter 2 and 10 of \cite{Wei94}. \par\indent
Let $\mathcal{A}$ be an abelian category with enough injectives, i.e. for every object $A\in \mathcal{A}$, there exists an injective object $I$ and a monomorphism $A\hookrightarrow I$. In particular, this implies that every object $A$ has an injective resolution $A\hookrightarrow I^{\bullet}$. As an example, the abelian category $\mathrm{Coh}X$ on a projective variety $X$ has enough injectives.
\begin{mydef}
Let $F: \mathcal{A}\rightarrow \mathcal{B}$ be a left exact functor between two abelian categories with $\mathcal{A}$ having enough injectives. We can define the \emph{right derived functors} $R^iF\,(i\geq 0)$ as follows. For each $A\in \mathcal{A}$, choose an injective resolution $A\hookrightarrow I^{\bullet}$ and define
$$R^iF(A)=H^i(F(I^{\bullet})).$$
\end{mydef}
Note that since $0\rightarrow F(A)\rightarrow F(I_0)\rightarrow F(I_1)$ is exact, we always have $R^0F\cong F$. It is an important fact that the definition of a right derived functor does not depend on the choice of the injective resolution. The proof of this is essentially the following: for any two injective resolutions of an object $A$, there exists a `lift' from one to the other that is unique up to chain homotopy. \par\indent
Let $(X,\mathcal{O}_X)$ be a ringed space and let $\mathcal{F}, \mathcal{G}$ be $\mathcal{O}_X$-modules. We denote by $\mathrm{Hom}(\mathcal{F},\mathcal{G})$ the homomorphism group of $\mathcal{O}_X$-modules, and $\SHom(\mathcal{F},\mathcal{G})$ the sheaf Hom construction. For fixed $\mathcal{F}$, $\mathrm{Hom}(\mathcal{F}, \cdot)$ is a left exact functor from $\mathrm{Mod}_{\mathcal{O}_X}$ to $\mathrm{Ab}$, and $\SHom(\mathcal{F}, \cdot)$ is a left exact functor from $\mathrm{Mod}_{\mathcal{O}_X}$ to $\mathrm{Mod}_{\mathcal{O}_X}$. Therefore, we may define their respective right derived functors as $\mathrm{Ext}^i(\mathcal{F},\cdot)$ and $\SExt^i(\mathcal{F},\cdot)$. As another example, the global section functor $\Gamma(X,\cdot)$ is left exact from $\mathrm{Mod}_{\mathcal{O}_X}$ to $\mathrm{Ab}$, and its right derived functors are just the sheaf cohomologies $H^i$. It is a famous theorem that when $X$ is a quasicompact and separated scheme and when $\mathcal{F}\in \mathrm{Coh}X$, the derived functor definition of sheaf cohomology agrees with \v{C}ech cohmology. \\\par\indent
Now we describe the construction of (bounded)derived category of an abelian category $\mathcal{A}$, which consists of three stages:\par\indent
1) We first consider the category of bounded cochain complexes $C^b(\mathcal{A})$, whose objects are cochain complexes $E^{\bullet}$ such that $H^i(E)=0$ for all but finitely many $i$, and morphisms are cochain maps.\par\indent
2) The homotopy category $K^b(\mathcal{A})$ is defined to have the same objects as $C^b(\mathcal{A})$, but two morphisms $f^{\bullet}, g^{\bullet}: E^{\bullet}\rightarrow F^{\bullet}$ are identified if they are homotopic, i.e. if there exists maps $h^i: E^i\rightarrow F^{i-1}$ such that $f_i-g_i=d\circ h^i-h^{i+1}\circ d$. \par\indent
3) Finally, the bounded derived category $D^b(\mathcal{A})$ is defined by `inverting' all quasi-isomorphisms, i.e. chain maps that induce isomorphisms on each cohomology group. Formally, this process is called \emph{localization} of a category, see \cite[Section 10.3]{Wei94}.\par\indent
In fact, $D^b(\mathcal{A})$ can be shown to be a triangulated category, equipped with the standard shift functor and whose distinguished triangles are given by the mapping cone construction.\\\par\indent
Derived categories and derived functors, as their names suggest, are closely related. One way to motivate the derived functor construction from a derived category perspective is the following. If $F: \mathcal{A}\rightarrow \mathcal{B}$ is a functor between two abelian categoies, then $F$ naturally extends to functors $C^b(F): C^b(\mathcal{A})\rightarrow C^b(\mathcal{B})$ and $K^b(F): K^b(\mathcal{A})\rightarrow K^b(\mathcal{B})$. The reason is that the relations defining chain complexes and homotopies are both functorial. In contrast, however, $F$ does not itself define a functor from $D^b(\mathcal{A})$ to $D^b(\mathcal{B})$ unless $F$ is exact. But is there a natural way to extend $F$ to derived categories? The answer to this question is exactly(no pun intended) derived functors. \par\indent
Before we give the construction, we need the following proposition about injective objects in an abelian category.
\begin{prop}
Let $\mathcal{A}$ be an abelian category. Then the following hold:\\
1) If $A^{\bullet}\in D^b(\mathcal{A})$ and $I^{\bullet}$ a bounded complex consists of injectives, then
$$\mathrm{Hom}_{D^b(\mathcal{A})}(A^{\bullet}, I^{\bullet})=\mathrm{Hom}_{K^b(\mathcal{A})}(A^{\bullet}, I^{\bullet}).$$
2) Suppose $\mathcal{A}$ has enough injectives, then
$$D^b(\mathcal{A})\cong K^b(\mathcal{I}),$$
where $K^b(\mathcal{I})$ is the full subcategory of $K^b(\mathcal{A})$ whose objects are bounded complexes consisting of injectives.
\end{prop}
For a proof of this proposition, see \cite[Section 10.4]{Wei94}. In 2), the equivalence of categories is given by sending a complex to the total complex of its Cartan-Eilenberg resolution. Given this proposition, we can define the derived functor of $F: \mathcal{A}\rightarrow \mathcal{B}$ as the composition
$$D^b(\mathcal{A})\cong K^b(\mathcal{I})\rightarrow K^b(\mathcal{B})\rightarrow D^b(\mathcal{B}),$$
where the middle map is $K^b(F)$, and the last map is simply passing to the localization. It can be easily verified that for $A^{\bullet}\in D^b(\mathcal{A})$ concentrated in degree $0$, then the definition of derived functor we just gave is the same as the `naive' definition given at the beginning of this section. \par\indent
Another useful application of Proposition 5.2 is the following lemma.
\begin{lemma}
For $A\in \mathcal{A}$, $A[i]$ denote the complex whose $(-i)$-th entry is $A$ and zero everywhere else. Then, for any $E,F\in \mathcal{A}$, we have
$$\mathrm{Hom}_{D^b(\mathcal{A})}(E,F[i])=
\begin{cases}
0 & \textrm{for}\;i<0\\
\mathrm{Ext}^i(E,F) &\textrm{for}\;i\geq 0
\end{cases}
$$
\end{lemma}

\subsection*{Serre duality}
\begin{mydef}
Let $S$ be a graded ring and $(\mathrm{Proj}S, \mathcal{O}_{\mathrm{Proj}S})$ be the usual proj construction. Define the \emph{Serre twsiting sheaf} $\mathcal{O}(1)$ as $\widehat{M}$ where $M_n=S_{n+1}$. Similarly, define $\mathcal{O}(d)$ as $\widehat{M}$ where $M_n=S_{n+d}$.
\end{mydef}
In the above definition, it is easy to verify that $\mathcal{O}(d)=\mathcal{O}(1)^{\otimes d}$. Let $\mathbb{P}^n_{\mathbb{C}}=\mathrm{Proj}\,\mathbb{C}[x_0,\cdots,x_n]$ denote the projective $n$-space over $\mathbb{C}$.
\begin{thm}
(i) $H^0((\mathbb{P}^n, \mathcal{O}_{\mathbb{P}^n}(m))$ is a vector space of dimension $\binom{n+m}{m}$ if $m\geq 0$.\\
(ii) $H^n((\mathbb{P}^n, \mathcal{O}_{\mathbb{P}^n}(m))$ is a vector space of dimension $\binom{-m-1}{-n-m-1}$ if $m\leq -n-1$.\\
(iii) $H^i((\mathbb{P}^n, \mathcal{O}_{\mathbb{P}^n}(m))=0$ otherwise. \\
(iv) The natural map
$$H^0(\mathbb{P}^n, \mathcal{O}_{\mathbb{P}^n}(m))\times H^n(\mathbb{P}^n, \mathcal{O}_{\mathbb{P}^n}(-m-n-1))\rightarrow H^n(\mathbb{P}^n, \mathcal{O}_{\mathbb{P}^n}(-n-1))$$
is a perfect pairing.
\end{thm}
\noindent\emph{Proof.} See \cite[Chapter III, Theorem 5.1]{Har77}.\\\par\indent
In fact, Theorem 5.2(iv) is a special case of a general principle called Serre duality. Before stating Serre duality, we first recall that for any ringed space $X$, $\mathrm{Hom}(\mathcal{O}_X,\cdot)$ and $\Gamma(X,\cdot)$ represent the same functor. As a result, their derived functor are also identical, i.e. $\mathrm{Ext}^i(\mathcal{O}_X,\cdot)=H^i(X,\cdot)$.\par\indent
Let $X=\mathbb{P}^n_{\mathbb{C}}$, and let $\omega_X=\bigwedge^n\Omega_{X/k}\cong \mathcal{O}(-n-1)$ be the canonical bundle. By Theorem 5.2, we have $H^n(X, \omega_X)\cong \mathbb{C}$.
\begin{thm}
Let $\mathcal{F}$ be a coherent sheaf over $X$. Then, the natural pairing
$$\mathrm{Ext}^i(\mathcal{F}, \omega_X)\times H^{n-i}(X,\mathcal{F})\rightarrow H^n(X,\omega_X)\cong \mathbb{C}$$
is perfect for all $0\leq i\leq n$. As a consequence, there is a natural isomorphism
$$\mathrm{Ext}^i(\mathcal{F}, \omega_X)\cong H^{n-i}(X,\mathcal{F})^*$$
for each $0\leq i\leq n$.
\end{thm}
\noindent\emph{Proof.} See \cite[Chapter III, Theorem 7.6]{Har77}.\\\par\indent
In the language of derived category, the above theorem says that $-\otimes \omega_X[n]$ is a Serre functor on $D^bCoh(X)$.

\end{document}